\documentclass[reqno,12pt]{article}
\usepackage{hyperref}
\usepackage{amsmath}
\usepackage{amssymb}
\usepackage{mathrsfs}

\newtheorem{thm}{Theorem}[section]
\newtheorem{lemma}[thm]{Lemma}

\newtheorem{corol}[thm]{Corollary}
\newtheorem{propos}[thm]{Proposition}
\newtheorem{rema}{Remark}[section]
\def\bp{\begin{propos}}
\def\ep{\end{propos}}
\def\bt{\begin{thm}}
\def\et{\end{thm}}
\def\bco{\begin{corol}}
\def\eco{\end{corol}}
\def\bl{\begin{lemma}}
\def\el{\end{lemma}}
\def\br{\begin{rema}}
\def\er{\end{rema}}
\def\be{\begin{equation}}
\def\ee{\end{equation}}
\def\ba{\begin{array}}
\def\ea{\end{array}}
\def\bena{\begin{eqnarray}}
\def\eena{\end{eqnarray}}

\def\P{{\mathbb P}}
\def\E{{\mathbb E}}
\def\R{{\mathbb R}}

\def\1{I}
\def\imath{\textbf{i}}
\def\jmath{\textbf{j}}
\def\chi{\zeta}

\def\QED{\hfill$\square$\vskip 3mm}

\def\Dp{\displaystyle}
\def\Df{\Dp\frac}

\def\({\left(}
\def\){\right)}

\begin{document}

\title{A note on the asymptotic behavior of the height for a birth-and-death process
\\[5mm]
\footnotetext{*Correspondence author}
\footnotetext{AMS classification (2010): 60F05, 60J27} 
\footnotetext{Research supported in part by the Natural Science Foundation of China (under
grants 11471222, 11671275) and Academy for Multidisciplinary Studies, Capital Normal University}}

\author{Feng Wang$^*$,\ \ \ Xian-Yuan Wu,\ \ \ Rui Zhu}\vskip 10mm
\date{}
\maketitle

\begin{center}
\begin{minipage}{13cm}

\noindent School of Mathematical Sciences, Capital
Normal University, Beijing, 100048, China. Email:
\texttt{wangf@mail.cnu.edu.cn},\ \texttt{wuxy@cnu.edu.cn},\ \texttt{1073755862@qq.com}
\end{minipage}
\end{center}
\vskip 5mm
{\begin{center} \begin{minipage}{13.5cm}
{\bf Abstract}: This paper focuses on the asymptotic behaviors of the {\it height} for
a birth-and-death process which related to a mean-field model \cite{FFS}(or the Anick-Mitra-Sondhi model \cite{DDM}). Recently,
the asymptotic mean value of the height for the model is given in \cite{LAV}. In this paper, first, the asymptotic variance of the height is given, and as a consequence, a weak Law of Large Number for the height is obtained. Second, the centered and normalized height is proved to converge in distribution to a degenerate distribution, this indicates that the desired Central Limit Theorem fails.

\end{minipage}
\end{center}}

 \vskip 5mm
\section{Introduction and statement of the results}
\renewcommand{\theequation}{1.\arabic{equation}}
\setcounter{equation}{0}

Birth-and-death process is a continuous-time Markov chain, which plays important roles in stochastic processes and queuing theory  \cite{RBC,LSC}.
Here we consider a birth-and-death process related to mean-field model and the Anick-Mitra-Sondhi model \cite{DDM}. Let $\{X_t,t\geq0\}$ be the birth-and-death process with state space $E=\{0,1,2,\cdots,N\}$ and the following conservative $Q$-matrix
$$    
Q=\left(                 
\begin{array}{cccccccc}   
q_{0,0} & N\nu & 0 & 0 & \cdots & 0 & 0 & 0 \\  
\mu & q_{1,1} & (N-1)\nu & 0 & \cdots & 0 & 0 & 0 \\  
0 & 2\mu & q_{2,2}& (N-2)\nu & \cdots & 0 & 0 & 0 \\
\vdots & \vdots & \vdots & \vdots &   & \vdots & \vdots & \vdots \\
 0 & 0 & 0 & 0 & \cdots & (N-1)\mu & q_{N-1,N-1}& \nu \\
  0 & 0 & 0 & 0 & \cdots & 0 & N\mu & q_{N,N}
\end{array}
\right),                 
$$
where $\mu,\nu>0$. Note that a $Q$-matrix is called conservative, if its row summation is zero.

Let $\rho=\dfrac{\nu}{\mu}$. Clearly, the chain $\{X_t:t\geq0\}$ is ergodic and has stationary distribution
\be
\pi_k:=\dfrac{1}{(1+\rho)^N}\binom{N}{k}\rho^k,\ k\in E.
\ee
Note that the transition probability matrix of its jump chain is given by
$$    
P=\left(                 
\begin{array}{cccccccc}   
0 & 1 & 0 & 0 & \cdots & 0 & 0 & 0 \\  
\frac{1}{1+(N-1)\rho} & 0 & \frac{(N-1)\rho}{1+(N-1)\rho}  & 0 & \cdots & 0 & 0 & 0 \\  
0 & \frac{2}{2+(N-2)\rho} & 0 & \frac{(N-2)\rho}{2+(N-2)\rho}  & \cdots & 0 & 0 & 0 \\
\vdots & \vdots & \vdots & \vdots &   & \vdots & \vdots & \vdots \\
0 & 0 & 0 & 0 & \cdots &\frac{N-1}{N-1+\rho} & 0 & \frac{\rho}{N-1+\rho} \\
0 & 0 & 0 & 0 & \cdots & 0 & 1 & 0
\end{array}
\right).                 
$$

The process $\{X_t,t\geq0\}$ has been studied in statistical physics
as a mean-field model (see \cite{FFS}), and as the Anick-Mitra-Sondhi model \cite{DDM} which is related to data-handling system with multiple sources. When $N=2$ and $N=3$, $\{X_t,t\geq0\}$ is also considered as a genomic model in
\cite{hsv}.

The present paper focuses on the {\it height} of $\{X_t:t\geq0\}$. Write $\{t:X_t>0\}=\cup_{i=1}^{\infty}[\tau_i,\eta_i)$, where $\{[\tau_i,\eta_i),i\in \mathbb{N}\}$ is the family of maximal disjoint time intervals such that $X_t>0$ on every interval $[\tau_i,\eta_i)$.
We consider the random variables
$$H_N^{(i)}:=\max\{X_t,t\in[\tau_i,\eta_i)\},$$
the possible values of $H_N^{(i)}$ may be listed as $1,2,\cdots,N$.
By definition of $H_N^{(i)}$, $X_{\tau_i}=1$,  $H_N^{(i)}$ is the maximal value
which $X_t$ can reach, before return to 0.  $H_N^{(i)}$ is called the height of $[\tau_i,\eta_i)$, and can be regarded as the maximal number of jobs that are served concurrently during a busy period in task-allocation problems,
or the maximal number of occupied nodes in a mean-field model before all nodes are free.
 $\{H_N^{(i)},i\in \mathbb{N}\}$  are
independent and identically distributed. This is due to the fact that
$X_t$ is Markov chain and $X_{\eta_i}=0$. The distribution of $H_N^{(i)}$
does not depend on $i$, $H_N^{(i)}$ is reduced to $H_N$. The asymptotic behavior of $H_N$  is studied in the case when the number of states tends to infinity.

The asymptotic mean value of $H_N$
is considered in \cite{LAV}:
%

\bt$^{\cite[\rm{Theorem\ 1}]{LAV}}$\label{thA}
For $ \rho\in (0,1)$, let $\alpha:=\alpha(\rho)$ be the unique solution of
the equation $x^{x}(1-x)^{1-x}=\rho^{x}$, let
\be\label{1.0}
f(\rho)=\left\{
\begin{array}{rcl}	
\alpha, & &0<\rho<1,\\
1, & & \rho\geq1.
\end{array} \right.
\ee
Then,	
\be\label{1.1}
	\lim_{N\rightarrow\infty}\frac{\mathbb{E}\left(H_N\right)}{N}=f(\rho).
\ee
\et

In the present paper, following the work of \cite{LAV}, we study and the fluctuations for $H_N$. Firstly, we have the following asymptotic behavior of the variance of $H_N$.

\bt\label{th1}
Suppose that $f(\rho)$ is given in Theorem\ref{thA}, then
\be\label{1.2}
	\lim_{N\rightarrow\infty}\frac{{\rm{Var}}\left(H_N\right)}{N}
	=\Df{f^2(\rho)}{\rho},
\ee and
\be\label{1.3}
	\lim_{N\rightarrow\infty}\frac{H_N}{N}={f(\rho)} \ \ \ \  \ in \ \ probability.	
\ee
\et
Secondly, we give a upper bound to the fluctuation of $H_N$ as follows.
\bt\label{th2}	Suppose $\varphi(N)$ satisfies that
	$\lim_{N\rightarrow\infty}\dfrac{\log N}{\varphi(N)}=0$, then
	\be\label{1.4}
	\lim_{N\rightarrow\infty}\mathbb{P}\left( \frac{H_N-\mathbb{E}\left(H_N\right)}{\varphi(N)}\leq  x\right)=\left\{
	\begin{array}{rcl}	
	0, & &x<0,\\
	1, & & x> 0.
	\end{array} \right.	
	\ee
\et	

\br
Theorem\ref{th2} indicates that the fluctuation for $H_N$ is upper bounded by $\log N$. In the case when $\varphi(N)=\sqrt{{\rm Var}\left(H_N\right)}$, (\ref{1.4}) shows that the centered and normalized $H_N$ converges weakly to a degenerate random variable.
\er

\vskip 5mm
\section{Proofs}
\renewcommand{\theequation}{2.\arabic{equation}}
\setcounter{equation}{0}
Before giving proofs, we give some useful notations. For any $x\in \R$, denote by $[x]$ the integer part of $x$. For positive series $\{a_n:n\geq 1\}$ and $\{b_n:n\geq 1\}$, write $b_n=O(a_n)$ if there exists some constant $C>0$ such that $b_n\leq Ca_n$ for all large enough $n$; write $b_n=\Theta (a_n)$, if $b_n=O(a_n)$ and $a_n=O(b_n)$.

By the law of total probability and iteration, the distribution of $H_N$ is given in \cite{LAV} as the following:

\bl$^{\cite[{\rm Lemma\ 1}]{LAV}}$\label{l1}

\be\label{2.0}	\P\(H_N\geq k\)=\Df{1}{\Dp\sum_{i=0}^{k-1}\Df{1}{\rho^i\binom{N-1}{i}}},k=1,2,\ldots,N.\ee
\el
	
Let $r_{\rho,n}(i):=\rho^{-i}\binom{n-1}{i}^{-1},$ $i=1,2,\ldots,n-1$. It is straightforward to check that $r_{\rho,n}(i)$ decreases strictly in $i$ when $i<\left[{n}/{(\rho+1)}\right]$, and increases strictly in $i$ otherwise.
For $0<\rho<1$, it was proved in \cite{LAV} that $\alpha$ is the unique solution of equation $x^x(1-x)^{1-x}=\rho^x$,
and $\rho<\alpha<1.$ Let $h_n=[\alpha(n-1)]$, then by Stirling's formula, we have
\be\label{2.1}
	r_{\rho,n}(h_n)=\Theta (\sqrt{n}).
\ee
Before giving proofs to the theorems, we shall give the following lemmas.

\bl\label{l2} Let $\rho\in (0,1)$, then for constants
	$C_1=\frac{2}{\log\alpha-\log\rho(1-\alpha)} $ and
	$ C_2=\frac{3}{\log\alpha-\log\rho}$, we have
	\begin{align}
		r_{\rho,n}\left(h_n+[C_1\log n]\right)&\geq r_{\rho,n}(h_n)n^2,\label{2.2}\\
		r_{\rho,n}\left(h_n-[C_2\log n]\right)&\leq r_{\rho,n}(h_n)n^{-3}.\label{2.3}
	\end{align}
\el
{\it Proof.} First, we prove (\ref{2.2}). By the definition of $r_{\rho,n}(i)$,
	\begin{align*}
		\frac{r_{\rho,n}(h_n+[C_1\log n])}{r_{\rho,n}(h_n)}&=
		\prod_{i=1}^{[C_1\log n]}\frac{r_{\rho,n}(h_n+i)}{r_{\rho,n}(h_n+i-1)}=\prod_{i=1}^{[C_1\log n]}\rho^{-1}\frac{h_n+i}{n-h_n-i-1}\\[2mm]
        &\geq\left(\frac{\alpha}{\rho(1-\alpha)}\right)^{[C_1\log n]}
        \cdot\prod_{i=1}^{[C_1\log n]}\dfrac{1+\frac{i-1}{\alpha n}}{1-\frac{i}{(1-\alpha)(n-1)}}\\[2mm]
        &\geq\left(\frac{\alpha}{\rho(1-\alpha)}\right)^{[C_1\log n]}=n^2.
   	\end{align*}
 Second, we obtain (\ref{2.3}) as the following:       	
\begin{align*}
	\frac{r_{\rho,n}(h_n-[C_2\log n])}{r_{\rho,n}(h_n)}&=
	\prod_{i=1}^{[C_2\log n]}\frac{r_{\rho,n}(h_n-i)}{r_{\rho,n}(h_n-i+1)}=\prod_{i=1}^{[C_2\log n]}\rho\frac{n-h_n+i-2}{h_n-i+1}\\[2mm]
	&\leq\left(\frac{\rho}{\alpha}\right)^{[C_2\log n]}
\prod_{i=1}^{[C_2\log n]}\frac{1-\alpha+\frac{i}{n-1}}{1-\frac{i-1}{\alpha (n-1)}}\\[2mm]
&\leq\left(\frac{\rho}{\alpha}\right)^{[C_2\log n]}\leq n^{-3}.
\end{align*}
\QED

\bl\label{l3}
	For $0<\rho<1$, $C_3=\frac{\alpha(3+\rho)}{\rho^2}$ and $C_2$ as given in Lemma~\ref{l2}, we have
	\be\label{2.01}[\alpha N]-[C_2\log N]-C_3\leq \mathbb{E}\left(H_N\right)\leq [\alpha N]+1\ee
	for $N$ large enough. For $\rho\geq 1$, we have
\be\label{2.02}N-4\leq \mathbb{E}\left(H_N\right)\leq N\ee for $N$ large enough.
\el

{\it Proof.}
	First, we prove the lower bound part of (\ref{2.01}). By (\ref{2.3}), $r_{\rho,n}(h_n-[C_2\log n])\leq r_{\rho,n}(h_n)n^{-3}$.
	For $2\leq i\leq h_n-[C_2\log n]$, we have $r_{\rho,n}(i)\leq\frac{2}{\rho^2(n-1)(n-2)}$, then
	\begin{align*}
	\sum_{i=0}^{h_n-[C_2\log n]}r_{\rho,n}(i)&\leq1+\frac{1}{\rho(n-1)}
	+\frac{2}{\rho^2(n-1)(n-2)}\cdot\left(h_n-[C_2\log n]\right)\\
	&\leq1+\frac{1}{\rho(n-1)}	+\frac{2}{\rho^2(n-2)}\\[2mm]
	&\leq\dfrac{\rho(n-1)+1+3/\rho}{\rho(n-1)}.
	\end{align*}
Hence
$$
\ba{cl}
	\mathbb{P}\left(H_N\geq h_N-[C_2\log N]\right)
	&=\dfrac{1}{\sum_{i=0}^{h_N-[C_2\log N]}r_{\rho,N}(i)}\\[4mm]
	&\geq \dfrac{\rho(N-1)}{\rho(N-1)+1+3/\rho}\\[3mm]
&\geq 1-\Df{3+\rho}{ (N-1)\rho^2}.
\ea
$$
Thus, we have
	\be\label{2.4}\ba{cl}
	\mathbb{E}\left(H_N\right)&=\Dp\sum_{i=1}^{N}\mathbb{P}\left(H_N\geq i\right)\geq\sum_{i=1}^{h_N-[C_2\log N]}\mathbb{P}\left(H_N\geq i\right)\\
	&\geq (h_N-[C_2\log N])\mathbb{P}\left(H_N\geq h_N-[C_2\log N]\right)\\[2mm]
	&\geq\left(h_N-[C_2\log N]\right)\left(1-\Df{3+\rho}{ (N-1)\rho^2}\right)\\
	&\geq h_N-[C_2\log N]-C_3.
\ea\ee

	Second, we prove the upper bound part of (\ref{2.01}). For $i\geq h_n$, then $i>\left[{(n-1)}/{(1+\rho)}\right]$. Noticing the fact that $r_{\rho,n}(i)$ strictly increases in $i$, we have $r_{\rho,n}(i)\geq r_{\rho,n}(h_n)$. Hence, for $k\geq1$, we have
	$$\sum_{i=0}^{h_n+k-1}r_{\rho,n}(i)\geq\sum_{i=h_n}^{h_n+k-1}r_{\rho,n}(i)
	\geq kr_{\rho,n}(h_n).$$
	So that
	\be\label{2.5}
		\mathbb{P}\left(H_N\geq h_N+k\right)\leq\frac{1}{kr_{\rho,N}(h_N)},
	\ee
and
	\begin{align*}
	\mathbb{E}\left(H_N\right)&=\sum_{i=1}^{N}\mathbb{P}\left(H_N\geq i\right)\\
	&\leq h_N+\sum_{i=h_N+1}^{N}\mathbb{P}\left(H_N\geq i\right)\\
	&\leq h_N+\sum_{k=1}^{N-h_N}\frac{1}{kr_{\rho,N}(h_N)}.
	\end{align*}
	
By the relation	between harmonic series and natural logarithm, we have
	 $$\lim_{N\rightarrow\infty}
	 \left[\sum_{k=1}^{N}\frac{1}{k}-\log N\right]=\gamma,$$
	 where $\gamma$ is  Euler-Mascheroni constant. 	
	By (\ref{2.1}), we have
	$$\sum_{k=1}^{N-h_N}\frac{1}{kr_{\rho,N}(h_N)}=O\left(\frac{\log N}{\sqrt{N}}\right), $$
	then,
	\be\label{2.6}\mathbb{E}\left(H_N\right)\leq h_N+1,\ee for $N$ large enough. The inequality (\ref{2.01}) follows from (\ref{2.4}) and (\ref{2.6}).

For $\rho\geq1$, note that $r_{\rho,n}(i)\leq\binom{n-1}{i}^{-1}$, then
 $$\sum_{i=0}^{n-2}r_{\rho,n}(i)\leq\sum_{i=0}^{n-2}\binom{n-1}{i}^{-1}\leq1+\frac{3}{n-1},$$
 and then
 \begin{align}\label{2.7}
 \mathbb{P}\left(H_N\geq N-1\right)=\dfrac{1}{\sum_{k=0}^{N-2}r_{\rho,N}(k)}\geq \frac{N-1}{N+2}.
 \end{align}
By (\ref{2.7}), we obtain (\ref{2.02}) and finish the proof of the lemma as the following.
$$N\geq
\E\left(H_N\right)=\sum_{i=1}^{N}
\P\left(H_N\geq i\right)
\geq (N-1)\P\left(H_N\geq N-1\right)\geq N-4.$$

\QED

{\it Proof of Theorem~\ref{th1}.} For $\rho\geq1$, first we have
\begin{align*}
{\rm Var}(H_N)&=\sum_{i=1}^{N}\left(i-\mathbb{E}\left(H_N\right)\right)^2\mathbb{P}\left(H_N=i\right)\\
&\geq\left(1-\mathbb{E}\left(H_N\right)\right)^2\mathbb{P}\left(H_N=1\right),
\end{align*}
then by (\ref{2.02})
       \be\label{2.03}{\rm Var}(H_N)\geq\frac{(N-3)^2}{1+\rho(N-1)}.\ee

Second, let $c=c(\rho)$ be the constant such that $r_{\rho,N}(i)N^2\leq\dfrac{c}{N}$ for all $3\leq i\leq N-4$. Using the fact that
\be\label{2.07}\mathbb{P}(H_N=i)=\mathbb{P}\left(H_N\geq i\right)-\mathbb{P}(H_N\geq i+1)\leq r_{\rho,N}(i),\ee
we have
\be\label{2.04}
\ba{cl}
{\rm Var}\left(H_N\right )&=\Dp\sum_{i=1}^{N}\left[i-\mathbb{E}\left(H_N\right)\right]^2\mathbb{P}\left(H_N=i\right)\\
            &\leq \Dp N^2r_{\rho,N}(1)+ N^2r_{\rho,N}(2)+\sum_{i=3}^{N-4}
             N^2r_{\rho,N}(i)+\sum_{i=N-3}^{N}(N-i)^2\\
            &\Dp\leq\frac{N}{\rho}+\frac{3c}{\rho^2}+13.
\ea
\ee
Then (\ref{1.2}) follows from (\ref{2.03}) and (\ref{2.04}).

For $0<\rho<1$, first, by the lower bound given in (\ref{2.01}) we have
\be\label{2.05}\ba{cl}{\rm Var}\left(H_N\right)&\geq\left[1-\mathbb{E}\left(H_N\right)\right]^2\mathbb{P}\left(H_N=1\right)\\[2mm]&\geq([\alpha N]-[C_2\log N]-C_3)^2\cdot\dfrac{1}{1+\rho N}.\ea\ee
Second, by (\ref{2.01}) and (\ref{2.07}), we have
\begin{align*}
&\sum_{i=1}^{h_N-C_2\log N}\left[i-\mathbb{E}\left(H_N\right)\right]^2\mathbb{P}\left(H_N=i\right)\\
&\leq \alpha^2N^2r_{\rho,N}(1)+ \alpha^2N^2r_{\rho,N}(2)+\sum_{i=3}^{h_N-C_2\log N}
\alpha^2N^2r_{\rho,N}(i)\\
&\leq\frac{\alpha^2N}{\rho}+\frac{3\alpha^2}{\rho^2}+\frac{13\alpha^2}{N},
\end{align*}

$$\sum_{i=h_N-C_2\log N}^{h_N+C_1\log N}\left[i-\mathbb{E}\left(H_N\right)\right]^2\mathbb{P}
\left(H_N=i\right)
\leq (C_1+C_2)^2(\log N)^2$$

and

\begin{align*}
\sum_{i=h_N+C_1\log N}^{N}\left[i-\mathbb{E}\left(H_N\right)\right]^2\mathbb{P}
\left(H_N=i\right)
&\leq N^2\sum_{i=h_N+C_1\log N}^{N}\mathbb{P}\left(H_N=i\right)\\
&=N^2\mathbb{P}(H_N\geq h_N+C_1\log N)\\
&\leq N^2\frac{1}{r_{\rho,N}(h_N+C_1\log N)}\\
&\leq O\(\frac{1}{\sqrt{N}}\).
\end{align*}
Note that last inequality follows from (\ref{2.1}) and (\ref{2.2}).
Thus
\be\label{2.06}
\ba{cl}
&\Dp{\rm Var}(H_N)=\Dp\sum_{i=1}^{N}\left[i-\mathbb{E}\left(H_N\right)\right]^2\mathbb{P}\left(H_N=i\right)\\[2mm]
&\leq\Dp\frac{\alpha^2N}{\rho}+\frac{3\alpha^2}{\rho^2}+\frac{13\alpha^2}{N}+(C_1+C_2)^2(\log N)^2+O\(\frac{1}{\sqrt{N}}\).
\ea
\ee
The equation (\ref{1.2}) follows from (\ref{2.05}) and (\ref{2.06}).

Finally, by using Chebyshev's inequality, we prove (\ref{1.3}). Actually, for any $\varepsilon>0$
\begin{align*}
	\mathbb{P}\left(\left|\frac{H_N}{N}-{f(\rho)}\right|\geq\varepsilon\right)\leq&\frac{1}{\varepsilon^2}
	\mathbb{E}\left[\left(\frac{H_N}{N}-{f(\rho)}\right)^2\right]\\
	\leq&\frac{1}{\varepsilon^2}\left[\frac{{\rm Var}(H_N)}{N^2}+\left(\frac{\mathbb{E}\left(H_N\right)}{N}-{f(\rho)}\right)^2\right].
\end{align*}
Then, by (\ref{1.1}) and (\ref{1.2}), we have
$$\lim_{N\rightarrow\infty}\left[\frac{{\rm Var}(H_N)}{N^2}+\left(\frac{\mathbb{E}\left(H_N\right)}{N}-{f(\rho)}\right)^2\right]=0.$$
Thus
$$\frac{H_N}{N}\rightarrow{f(\rho)}\ \ {\rm as}\ N\rightarrow\infty$$ in probability.
\QED

{\it Proof of Theorem~\ref{th2}.}
 By the condition that
 $\lim_{N\rightarrow\infty}{\log N}/{\varphi(N)}=0$ and Lemma \ref{l3},  for any $ x>0$ and $N$ large enough, we have $x\varphi(N)\geq[C_2\log N]+C_3$, then
$$\ba{cl}\mathbb{P}\left(H_N\leq\mathbb{E}\left(H_N\right)+x\varphi(N)\right)&\geq\mathbb{P}\left(H_N\leq h_N-[C_2\log N]-C_3+x\varphi(N)\right)\\[2mm]&\geq\mathbb{P}\left(H_N\leq h_N\right).\ea$$
By (\ref{2.5}), we have
$$\lim_{N\rightarrow\infty}\mathbb{P}\left(H_N\leq h_N\right)=1$$
then
$$\lim_{N\rightarrow\infty}\mathbb{P}
\left(\dfrac{H_N-\mathbb{E}\left(H_N\right)}{\varphi(N)}\leq x\right)=1.$$

 For any $x<0$, for $N$ large  enough,  $x\varphi(N)\leq-[C_2\log N]-1$, then
$$\ba{cl}\mathbb{P}\left(H_N\leq\mathbb{E}\left(H_N\right)+x\varphi(N)\right)&\leq\mathbb{P}\left(H_N\leq h_N+x\varphi(N)\right)\\[2mm]&\leq\mathbb{P}\left(H_N\leq h_N-[C_2\log N]-1\right).\ea$$
By (\ref{2.4}), we have
$$\lim_{N\rightarrow\infty}\mathbb{P}\left(H_N\leq  h_N-C_2\log N\right)=0,$$
then
$$\lim_{N\rightarrow\infty}\mathbb{P}
\left(\dfrac{H_N-\mathbb{E}\left(H_N\right)}{\varphi(N)}\leq x\right)=0.$$
\QED


\vskip10mm

\end{document}